\documentclass[10pt,draft]{article}
\usepackage{amsfonts,amsmath,amssymb,cite}
\newcommand{\Imag}{\mathop{\textrm{Im}}\nolimits}
\newcommand{\Li}{\mathop{\textrm{Li}}\nolimits}
\renewcommand{\theequation}{\thesection.\arabic{equation}}
\begin{document}
\title{On the derivatives
$\partial^{2}P_{\nu}(z)/\partial\nu^{2}$ and $\partial
Q_{\nu}(z)/\partial\nu$\\ 
of the Legendre functions with respect to their degrees}
\author{Rados{\l}aw Szmytkowski \\*[3ex]
Atomic and Optical Physics Division, \\
Department of Atomic, Molecular and Optical Physics, \\
Faculty of Applied Physics and Mathematics,
Gda{\'n}sk University of Technology, \\
ul.\ Gabriela Narutowicza 11/12, 
80--233 Gda{\'n}sk, Poland \\
email: radoslaw.szmytkowski@pg.edu.pl}
\date{}
\maketitle
\begin{center}
\textbf{Published as: Integral Transforms Spec.\ Funct. 28 (2017)
645--62} \\*[1ex]
\textbf{doi: 10.1080/10652469.2017.1339039} \\*[5ex]
\end{center}
\begin{abstract} 
We provide closed-form expressions for the degree-derivatives
$[\partial^{2}P_{\nu}(z)/\partial\nu^{2}]_{\nu=n}$ and $[\partial
Q_{\nu}(z)/\partial\nu]_{\nu=n}$, with $z\in\mathbb{C}$ and
$n\in\mathbb{N}_{0}$, where $P_{\nu}(z)$ and $Q_{\nu}(z)$ are the
Legendre functions of the first and the second kind, respectively.
For $[\partial^{2}P_{\nu}(z)/\partial\nu^{2}]_{\nu=n}$, we find that
\begin{displaymath}
\frac{\partial^{2}P_{\nu}(z)}{\partial\nu^{2}}\bigg|_{\nu=n}
=-2P_{n}(z)\Li_{2}\frac{1-z}{2}+B_{n}(z)\ln\frac{z+1}{2}+C_{n}(z),
\end{displaymath}
where $\Li_{2}[(1-z)/2]$ is the dilogarithm function, $P_{n}(z)$ is
the Legendre polynomial, while $B_{n}(z)$ and $C_{n}(z)$ are certain
polynomials in $z$ of degree $n$. For $[\partial
Q_{\nu}(z)/\partial\nu]_{\nu=n}$ and $z\in\mathbb{C}\setminus[-1,1]$,
we derive 
\begin{eqnarray*}
\frac{\partial Q_{\nu}(z)}{\partial\nu}\bigg|_{\nu=n} 
&=& -P_{n}(z)\Li_{2}\frac{1-z}{2}
-\frac{1}{2}P_{n}(z)\ln\frac{z+1}{2}\ln\frac{z-1}{2}
+\frac{1}{4}B_{n}(z)\ln\frac{z+1}{2} 
\nonumber \\
&& -\,\frac{(-1)^{n}}{4}B_{n}(-z)\ln\frac{z-1}{2}
-\frac{\pi^{2}}{6}P_{n}(z)
+\frac{1}{4}C_{n}(z)-\frac{(-1)^{n}}{4}C_{n}(-z).
\end{eqnarray*}
A counterpart expression for $[\partial
Q_{\nu}(x)/\partial\nu]_{\nu=n}$, applicable when $x\in(-1,1)$, is
also presented. Explicit representations of the polynomials
$B_{n}(z)$ and $C_{n}(z)$ as linear combinations of the Legendre
polynomials are given.
\end{abstract}
%
%
\section{Introduction}
\label{I}
\setcounter{equation}{0}
Over the past 20 years or so, a growth of interest in parameter
derivatives of various special functions has been observed. The
research done on the subject is documented in a number of papers
reporting diverse methods for finding such derivatives for orthogonal
polynomials in one \cite{Froe94,Koep97,Szmy09a} and two
\cite{Akta14,Akta15,Akta16} variables, for Bessel functions
\cite{Bryc05,Sesm14,Duns15,Bryc16}, for Legendre and allied functions
\cite{Szmy06,Szmy07,Szmy09b,Bryc10,Cohl10,Cohl11,Szmy11,Szmy12}, and
also for various types of hypergeometric functions
\cite{Froe94,Abad03,Anca08,Anca09,Anca10}.

In Refs.\ \cite{Szmy06,Szmy07}, we presented results of our
investigations on the first-order derivative of the Legendre function
of the first kind with respect to its degree. We showed that
$[\partial P_{\nu}(z)/\partial\nu]_{\nu=n}$, with
$n\in\mathbb{N}_{0}$, is of the form
\begin{equation}
\frac{\partial P_{\nu}(z)}{\partial\nu}\bigg|_{\nu=n}
=P_{n}(z)\ln\frac{z+1}{2}+R_{n}(z),
\label{1.1}
\end{equation}
where $P_{n}(z)$ is the Legendre polynomial of degree $n$ and
$R_{n}(z)$ is another polynomial in $z$ of the same degree. We
investigated properties of the polynomials $R_{n}(z)$ and arrived at
their several explicit representations, including the following one:
\begin{equation} 
R_{n}(z)=2[\psi(2n+1)-\psi(n+1)]P_{n}(z)
+2\sum_{k=0}^{n-1}(-1)^{n+k}\frac{2k+1}{(n-k)(n+k+1)}P_{k}(z),
\label{1.2} 
\end{equation}
where $\psi(z)=\mathrm{d}\ln\Gamma(z)/\mathrm{d}z$ is the digamma
function.

In the year 2012, Dr.\ George P.\ Schramkowski kindly informed the
present author that in the course of doing research on a certain
problem in theoretical hydrodynamics, he had come across higher-order
derivatives $[\partial^{k}P_{\nu}(z)/\partial\nu^{k}]_{\nu=n}$, with
$n\in\mathbb{N}_{0}$ and $k\geqslant2$. Using \emph{Mathematica\/},
Schramkowski found that
\begin{equation}
\frac{\partial^{2}P_{\nu}(z)}{\partial\nu^{2}}\bigg|_{\nu=0}
=-2\Li_{2}\frac{1-z}{2},
\label{1.3}
\end{equation}
where 
\begin{equation}
\Li_{2}z=-\int_{0}^{z}\mathrm{d}t\:\frac{\ln(1-t)}{t}
\label{1.4}
\end{equation}
is the dilogarithm function \cite{Lewi81,Apos10}. In Ref.\
\cite{Szmy13}, we gave an analytical proof of the result displayed in
Eq.\ (\ref{1.3}), and also we derived a closed-form formula for the
third-order derivative
$[\partial^{3}P_{\nu}(z)/\partial\nu^{3}]_{\nu=0}$. That work was
then extended by Laurenzi \cite{Laur15}, who found an expression for
the fourth-order derivative
$[\partial^{4}P_{\nu}(z)/\partial\nu^{4}]_{\nu=0}$.

The primary purpose of the present work is to pursue further the
research initiated by Schramkowski and continued by us
in Ref.\ \cite{Szmy13}. We shall show that for arbitrary
$n\in\mathbb{N}_{0}$ the second-order derivative
$[\partial^{2}P_{\nu}(z)/\partial\nu^{2}]_{\nu=n}$ may be expressed
in the form
\begin{equation}
\frac{\partial^{2}P_{\nu}(z)}{\partial\nu^{2}}\bigg|_{\nu=n}
=-2P_{n}(z)\Li_{2}\frac{1-z}{2}+B_{n}(z)\ln\frac{z+1}{2}+C_{n}(z),
\label{1.5}
\end{equation}
where the polynomials $B_{n}(z)$ and $C_{n}(z)$ have the following
representations in terms of the Legendre polynomials:
\begin{equation}
B_{n}(z)=4[\psi(2n+1)-\psi(n+1)]P_{n}(z)
+4\sum_{k=0}^{n-1}\frac{2k+1}{(n-k)(n+k+1)}P_{k}(z)
\label{1.6}
\end{equation}
and
\begin{eqnarray}
C_{n}(z) &=& \left\{-\frac{\pi^{2}}{3}+4[\psi(2n+1)-\psi(n+1)]^{2}
+4\psi_{1}(2n+1)-2\psi_{1}(n+1)\right\}P_{n}(z)
\nonumber \\
&& +\,4\sum_{k=0}^{n-1}(-1)^{n+k}\frac{2k+1}{(n-k)(n+k+1)}
\bigg\{2\bigg[\psi(n+k+1)-\psi(n-k+1)
\nonumber \\
&& \qquad -\,\psi\left(\left\lfloor\frac{n+k}{2}\right\rfloor+1\right)
+\psi\left(\left\lfloor\frac{n-k}{2}\right\rfloor+1\right)\bigg]
\nonumber \\
&& \quad -\,(-1)^{n+k}\frac{2k+1}{(n-k)(n+k+1)}
-\frac{2n+1}{(n-k)(n+k+1)}\bigg\}P_{k}(z),
\label{1.7}
\end{eqnarray}
with $\psi(z)$ being already defined under Eq.\ (\ref{1.2}),
$\psi_{1}(z)=\mathrm{d}\psi(z)/\mathrm{d}z$, and $\lfloor
x\rfloor=\max\{n\in\mathbb{Z}:\:n\leqslant x\}$ standing for the
integer part of $x$. In fact, the above result for
$[\partial^{2}P_{\nu}(z)/\partial\nu^{2}]_{\nu=n}$, valid for
$n\in\mathbb{N}_{0}$, may be easily extended to any $n\in\mathbb{Z}$,
since with the use of the well-known identity
\begin{equation}
P_{-\nu-1}(z)=P_{\nu}(z)
\label{1.8}
\end{equation}
one immediately finds that
\begin{equation}
\frac{\partial^{2}P_{\nu}(z)}{\partial\nu^{2}}\bigg|_{\nu=-n-1}
=\frac{\partial^{2}P_{\nu}(z)}{\partial\nu^{2}}\bigg|_{\nu=n}
\qquad (n\in\mathbb{N}_{0}).
\label{1.9}
\end{equation}

In addition to the above summarized study on the second-order
degree-derivative of the Legendre function of the first kind, which
will be presented in detail in Sec.\ \ref{II} below, later in Sec.\
\ref{III} we shall also prove that if $z\in\mathbb{C}\setminus[-1,1]$
and $n\in\mathbb{N}_{0}$, then the first-order derivative $[\partial
Q_{\nu}(z)/\partial\nu]_{\nu=n}$, where $Q_{\nu}(z)$ is the Legendre
function of the second kind, is given by
\begin{eqnarray}
\frac{\partial Q_{\nu}(z)}{\partial\nu}\bigg|_{\nu=n}
&=& -P_{n}(z)\Li_{2}\frac{1-z}{2}
-\frac{1}{2}P_{n}(z)\ln\frac{z+1}{2}\ln\frac{z-1}{2}
+\frac{1}{4}B_{n}(z)\ln\frac{z+1}{2}
\nonumber \\
&& -\,\frac{(-1)^{n}}{4}B_{n}(-z)\ln\frac{z-1}{2}
-\frac{\pi^{2}}{6}P_{n}(z)
+\frac{1}{4}C_{n}(z)-\frac{(-1)^{n}}{4}C_{n}(-z).
\label{1.10}
\end{eqnarray}
A counterpart expression for $[\partial
Q_{\nu}(x)/\partial\nu]_{\nu=n}$, applicable when $x\in(-1,1)$, will
also be derived.
%
%
\section{The derivatives
$[\partial^{2}P_{\nu}(z)/\partial\nu^{2}]_{\nu=n}$}
\label{II}
\setcounter{equation}{0}
\subsection{The general form of 
$[\partial^{2}P_{\nu}(z)/\partial\nu^{2}]_{\nu=n}$}
\label{II.1}
Our point of departure is the well-known recurrence relation
\begin{equation}
(\nu+1)P_{\nu+1}(z)-(2\nu+1)zP_{\nu}(z)+\nu P_{\nu-1}(z)=0
\label{2.1}
\end{equation}
obeyed by the Legendre functions of the first kind. Differentiating
it twice with respect to $\nu$ and setting then $\nu=n$ yields
\begin{eqnarray}
&& (n+1)\frac{\partial^{2}P_{\nu}(z)}
{\partial\nu^{2}}\bigg|_{\nu=n+1}
-(2n+1)z\frac{\partial^{2}P_{\nu}(z)}
{\partial\nu^{2}}\bigg|_{\nu=n}
+n\frac{\partial^{2}P_{\nu}(z)}
{\partial\nu^{2}}\bigg|_{\nu=n-1}
\nonumber \\
&& \qquad =\,-2\left[\frac{\partial P_{\nu}(z)}
{\partial\nu}\bigg|_{\nu=n+1}
-2z\frac{\partial P_{\nu}(z)}
{\partial\nu}\bigg|_{\nu=n}
+\frac{\partial P_{\nu}(z)}
{\partial\nu}\bigg|_{\nu=n-1}\right].
\label{2.2}
\end{eqnarray}
If we replace the first-order derivatives on the right-hand side with
expressions following from Eq.\ (\ref{1.1}), this furnishes
\begin{eqnarray}
&& (n+1)\frac{\partial^{2}P_{\nu}(z)}
{\partial\nu^{2}}\bigg|_{\nu=n+1}
-(2n+1)z\frac{\partial^{2}P_{\nu}(z)}
{\partial\nu^{2}}\bigg|_{\nu=n}
+n\frac{\partial^{2}P_{\nu}(z)}
{\partial\nu^{2}}\bigg|_{\nu=n-1}
\nonumber \\
&& \qquad =\,-2[P_{n+1}(z)-2zP_{n}(z)+P_{n-1}(z)]\ln\frac{z+1}{2}
-2[R_{n+1}(z)-2zR_{n}(z)+R_{n-1}(z)].
\nonumber \\
&&
\label{2.3}
\end{eqnarray}
From the formal point of view, one may look at Eq.\ (\ref{2.3}) as a
second-order difference equation and then two additional conditions
are necessary to single out the sequence
$[\partial^{2}P_{\nu}(z)/\partial\nu^{2}]_{\nu=n}$ from its general
solution. Such conditions may be chosen in a variety of ways but for
our purposes it is most convenient to take the explicit expression
for $[\partial^{2}P_{\nu}(z)/\partial\nu^{2}]_{\nu=0}$, given in Eq.\
(\ref{1.1}), as the first one. The second suitable condition follows
from Eq.\ (\ref{2.3}) after one lets $n=0$. Then, with the use of the
identities
\begin{equation}
P_{-1}(z)=1, 
\qquad P_{0}(z)=1,
\qquad P_{1}(z)=z
\label{2.4}
\end{equation}
and \cite[Sec.\ 5.2]{Szmy06}
\begin{equation}
R_{-1}(z)=-2\ln\frac{z+1}{2},
\qquad R_{0}(z)=0,
\qquad R_{1}(z)=z-1,
\label{2.5}
\end{equation}
Eq.\ (\ref{2.3}) reduces to the form
\begin{equation}
\frac{\partial^{2}P_{\nu}(z)}{\partial\nu^{2}}\bigg|_{\nu=1}
-z\frac{\partial^{2}P_{\nu}(z)}{\partial\nu^{2}}\bigg|_{\nu=0}
=2(z+1)\ln\frac{z+1}{2}-2(z-1).
\label{2.6}
\end{equation}
On combining Eq.\ (\ref{2.6}) with Eq.\ (\ref{1.1}), one finds that
\begin{equation}
\frac{\partial^{2}P_{\nu}(z)}{\partial\nu^{2}}\bigg|_{\nu=1}
=-2z\Li_{2}\frac{1-z}{2}+2(z+1)\ln\frac{z+1}{2}-2(z-1).
\label{2.7}
\end{equation}
If necessary, Eq.\ (\ref{2.3}) may be applied recursively, with Eqs.\
(\ref{1.3}) and (\ref{2.7}) used as initial conditions, to generate
the derivative in question for any particular $n\geqslant2$. However,
as we shall show below, it is also possible to obtain a closed-form
representation for
$[\partial^{2}P_{\nu}(z)/\partial\nu^{2}]_{\nu=n}$.

As the first step towards that goal, we observe that the structure of
Eq.\ (\ref{2.3}), together with explicit expressions for the
derivatives $[\partial^{2}P_{\nu}(z)/\partial\nu^{2}]_{\nu=0}$ and
$[\partial^{2}P_{\nu}(z)/\partial\nu^{2}]_{\nu=1}$ displayed in Eqs.\
(\ref{1.3}) and (\ref{2.7}), fix the form of
$[\partial^{2}P_{\nu}(z)/\partial\nu^{2}]_{\nu=n}$ to be
\begin{equation}
\frac{\partial^{2}P_{\nu}(z)}{\partial\nu^{2}}\bigg|_{\nu=n}
=A_{n}(z)\Li_{2}\frac{1-z}{2}+B_{n}(z)\ln\frac{z+1}{2}+C_{n}(z),
\label{2.8}
\end{equation}
where $A_{n}(z)$, $B_{n}(z)$ and $C_{n}(z)$ are polynomials in $z$ of
degree $n$. Since the right-hand side of Eq.\ (\ref{2.3}) does not
contain the dilogarithm function, the polynomial $A_{n}(z)$ solves
the homogeneous recurrence
\begin{equation}
(n+1)A_{n+1}(z)-(2n+1)zA_{n}(z)+nA_{n}(z)=0
\label{2.9}
\end{equation}
subject to the initial conditions
\begin{equation}
A_{0}(z)=-2=-2P_{0}(z),
\qquad
A_{1}(z)=-2z=-2P_{1}(z),
\label{2.10}
\end{equation}
which follow from Eqs.\ (\ref{1.1}) and (\ref{2.8}). Hence, we deduce
the following expression for $A_{n}(z)$ in terms of the Legendre
polynomial $P_{n}(z)$:
\begin{equation}
A_{n}(z)=-2P_{n}(z).
\label{2.11}
\end{equation}
Consequently, Eq.\ (\ref{2.8}) becomes
\begin{equation}
\frac{\partial^{2}P_{\nu}(z)}{\partial\nu^{2}}\bigg|_{\nu=n}
=-2P_{n}(z)\Li_{2}\frac{1-z}{2}+B_{n}(z)\ln\frac{z+1}{2}+C_{n}(z).
\label{2.12}
\end{equation}
The representations of the polynomials $B_{n}(z)$ and $C_{n}(z)$
remain to be established.
%
%
\subsection{Differential equations for the polynomials $B_{n}(z)$ and
$C_{n}(z)$}
\label{II.2}
It is known that the Legendre function $P_{\nu}(z)$ obeys the
differential identity
\begin{equation}
\left[\frac{\mathrm{d}}{\mathrm{d}z}(1-z^{2})
\frac{\mathrm{d}}{\mathrm{d}z}+\nu(\nu+1)\right]P_{\nu}(z)=0.
\label{2.13}
\end{equation}
If we differentiate it twice with respect to $\nu$ and then put
$\nu=n$, this gives
\begin{equation}
\left[\frac{\mathrm{d}}{\mathrm{d}z}(1-z^{2})
\frac{\mathrm{d}}{\mathrm{d}z}+n(n+1)\right]
\frac{\partial^{2}P_{\nu}(z)}{\partial\nu^{2}}\bigg|_{\nu=n}
=-2(2n+1)\frac{\partial P_{\nu}(z)}{\partial\nu}\bigg|_{\nu=n}
-2P_{n}(z)
\label{2.14}
\end{equation}
and further, after Eq.\ (\ref{1.1}) is plugged into the first term on
the right-hand side,
\begin{equation}
\left[\frac{\mathrm{d}}{\mathrm{d}z}(1-z^{2})
\frac{\mathrm{d}}{\mathrm{d}z}+n(n+1)\right]
\frac{\partial^{2}P_{\nu}(z)}{\partial\nu^{2}}\bigg|_{\nu=n}
=-2(2n+1)P_{n}(z)\ln\frac{z+1}{2}-2(2n+1)R_{n}(z)-2P_{n}(z).
\label{2.15}
\end{equation}
Next, we insert Eq.\ (\ref{2.12}) into the left-hand side of Eq.\
(\ref{2.15}) and equate those terms appearing on both sides which
involve the logarithmic factor. This yields the following
inhomogeneous differential equation for $B_{n}(z)$:
\begin{equation}
\left[\frac{\mathrm{d}}{\mathrm{d}z}(1-z^{2})
\frac{\mathrm{d}}{\mathrm{d}z}+n(n+1)\right]B_{n}(z)
=4\left[(z+1)\frac{\mathrm{d}P_{n}(z)}{\mathrm{d}z}-nP_{n}(z)\right].
\label{2.16}
\end{equation}
Similarly, after equating polynomial expressions on both sides, we
arrive at the inhomogeneous equation for $C_{n}(z)$:
\begin{equation}
\left[\frac{\mathrm{d}}{\mathrm{d}z}(1-z^{2})
\frac{\mathrm{d}}{\mathrm{d}z}+n(n+1)\right]C_{n}(z)
=2(z-1)\frac{\mathrm{d}B_{n}(z)}{\mathrm{d}z}+B_{n}(z)
-2(2n+1)R_{n}(z).
\label{2.17}
\end{equation}

Consider Eq.\ (\ref{2.16}). It is evident that it does not possess a
unique polynomial solution, since to any particular solution of that
form one may add an arbitrary multiple of the Legendre polynomial
$P_{n}(z)$, which results in another polynomial solution. To make the
polynomial solution unique, we thus need an additional constraint.
The latter follows from the limiting relation \cite{Magn66}
\begin{equation}
P_{\nu}(z)\stackrel{z\to-1}{\longrightarrow}
\frac{\sin(\pi\nu)}{\pi}\ln\frac{z+1}{2}+O(1),
\label{2.18}
\end{equation}
from which we find
\begin{equation}
\frac{\partial^{2}P_{\nu}(z)}{\partial\nu^{2}}\bigg|_{\nu=n}
\stackrel{z\to-1}{\longrightarrow}O(1).
\label{2.19}
\end{equation}
Hence, the left-hand side of Eq.\ (\ref{2.12}) remains finite for
$z\to-1$ and to make the right-hand side also finite in that limit,
we are forced to put
\begin{equation}
B_{n}(-1)=0.
\label{2.20}
\end{equation}
If, in turn, we wish to make the polynomial solution to Eq.\
(\ref{2.17}) unique, we use the identity
\begin{equation}
P_{\nu}(1)=1.
\label{2.21}
\end{equation}
Differentiating twice with respect to $\nu$, we obtain
\begin{equation}
\frac{\partial^{2}P_{\nu}(1)}{\partial\nu^{2}}\bigg|_{\nu=n}=0.
\label{2.22} 
\end{equation}
Since $\Li_{2}0=0$ and $\ln1=0$, we deduce that $C_{n}(z)$ is
constrained to obey
\begin{equation}
C_{n}(1)=0.
\label{2.23}
\end{equation}
Below we shall exploit Eqs.\ (\ref{2.16}), (\ref{2.20}), (\ref{2.17})
and (\ref{2.23}) to determine the polynomials $B_{n}(z)$ and
$C_{n}(z)$.
%
%
\subsection{Construction of the polynomials $B_{n}(z)$}
\label{II.3}
A general form of the polynomials $B_{n}(z)$ may be obtained with
ease. In Ref.\ \cite[Sec.\ 5.2.2]{Szmy06}, we have found that the
polynomials $R_{n}(z)$ obey the differential relation
\begin{equation}
\left[\frac{\mathrm{d}}{\mathrm{d}z}(1-z^{2})
\frac{\mathrm{d}}{\mathrm{d}z}+n(n+1)\right]R_{n}(z)
=2\left[(z-1)\frac{\mathrm{d}P_{n}(z)}{\mathrm{d}z}-nP_{n}(z)\right].
\label{2.24}
\end{equation}
Hence, with the use of the well-known identity
\begin{equation}
P_{n}(-z)=(-1)^{n}P_{n}(z),
\label{2.25}
\end{equation}
we deduce that
\begin{equation}
\left[\frac{\mathrm{d}}{\mathrm{d}z}(1-z^{2})
\frac{\mathrm{d}}{\mathrm{d}z}+n(n+1)\right]R_{n}(-z)
=(-1)^{n}2\left[(z+1)\frac{\mathrm{d}P_{n}(z)}{\mathrm{d}z}
-nP_{n}(z)\right].
\label{2.26}
\end{equation}
Comparison of Eqs.\ (\ref{2.16}) and (\ref{2.26}) shows that the
polynomial $B_{n}(z)$ must be of the form
\begin{equation}
B_{n}(z)=(-1)^{n}2R_{n}(-z)+b_{n}P_{n}(z).
\label{2.27}
\end{equation}
To determine the constant $b_{n}$, we put $z=-1$ in Eq.\
(\ref{2.27}). By virtue of the constraint (\ref{2.20}), with the use
of the relations
\begin{equation}
P_{n}(-1)=(-1)^{n}
\label{2.28}
\end{equation}
and (cf.\ Ref.\ \cite[Eq.\ (5.10)]{Szmy06})
\begin{equation}
R_{n}(1)=0,
\label{2.29}
\end{equation}
we infer that
\begin{equation}
b_{n}=0,
\label{2.30}
\end{equation}
and thus finally we arrive at
\begin{equation}
B_{n}(z)=(-1)^{n}2R_{n}(-z).
\label{2.31}
\end{equation}
On combining Eq.\ (\ref{2.31}) with Eq.\ (\ref{1.2}), we have the
following explicit representation of $B_{n}(z)$:
\begin{equation}
B_{n}(z)=4[\psi(2n+1)-\psi(n+1)]P_{n}(z)
+4\sum_{k=0}^{n-1}\frac{2k+1}{(n-k)(n+k+1)}P_{k}(z).
\label{2.32}
\end{equation}
Further expressions for $B_{n}(z)$ may be obtained if one combines
Eq.\ (\ref{2.31}) with Eqs.\ (2.1), (2.2), (2.4) and (5.90) from
Ref.\ \cite{Szmy06} or with Eqs.\ (11) and (12) from Ref.\
\cite{Szmy07}.
%
%
\subsection{Construction of the polynomials $C_{n}(z)$}
\label{II.4}
We shall seek a representation of $C_{n}(z)$ in the form of a linear
combination of Legendre polynomials:
\begin{equation}
C_{n}(z)=\sum_{k=0}^{n}c_{nk}P_{k}(z).
\label{2.33}
\end{equation}
Action on both sides of Eq.\ (\ref{2.33}) with the Legendre
differential operator appearing on the left-hand side of Eq.\
(\ref{2.17}) gives
\begin{equation}
\left[\frac{\mathrm{d}}{\mathrm{d}z}(1-z^{2})
\frac{\mathrm{d}}{\mathrm{d}z}+n(n+1)\right]C_{n}(z)
=\sum_{k=0}^{n-1}(n-k)(n+k+1)c_{nk}P_{k}(z).
\label{2.34}
\end{equation}
On the other side, with the aid of Eqs.\ (\ref{1.2}) and
(\ref{2.32}), and of the identity
\begin{equation}
(z-1)\frac{\mathrm{d}P_{n}(z)}{\mathrm{d}z}
=nP_{n}(z)+\sum_{k=0}^{n-1}(-1)^{n+k}(2k+1)P_{k}(z),
\label{2.35}
\end{equation}
after some algebra we find that the expression on the right-hand side
of Eq.\ (\ref{2.17}) may be written as
\begin{eqnarray}
&& \hspace*{-3em}
2(z-1)\frac{\mathrm{d}B_{n}(z)}{\mathrm{d}z}+B_{n}(z)-2(2n+1)R_{n}(z)
\nonumber \\
&=& \sum_{k=0}^{n-1}\bigg\{(-1)^{n+k}8(2k+1)[\psi(2n+1)-\psi(n+1)]
\nonumber \\
&& +\,8(2k+1)\sum_{m=k}^{n-1}(-1)^{m+k}\frac{2m+1}{(n-m)(n+m+1)}
\nonumber \\
&& -\,\frac{4(2k+1)^{2}}{(n-k)(n+k+1)}
-(-1)^{n+k}\frac{4(2n+1)(2k+1)}{(n-k)(n+k+1)}\bigg\}P_{k}(z).
\label{2.36}
\end{eqnarray}
Equations (\ref{2.17}), (\ref{2.34}) and (\ref{2.36}) yield
\begin{eqnarray}
c_{nk} &=& (-1)^{n+k}[\psi(2n+1)-\psi(n+1)]
\frac{8(2k+1)}{(n-k)(n+k+1)}
\nonumber \\
&& +\,\frac{8(2k+1)}{(n-k)(n+k+1)}
\sum_{m=k}^{n-1}(-1)^{m+k}\frac{2m+1}{(n-m)(n+m+1)}
\nonumber \\
&& -\,\frac{4(2k+1)^{2}}{(n-k)^{2}(n+k+1)^{2}}
-(-1)^{n+k}\frac{4(2n+1)(2k+1)}{(n-k)^{2}(n+k+1)^{2}}
\qquad (0\leqslant k\leqslant n-1).
\nonumber \\
&&
\label{2.37}
\end{eqnarray}
It is proven in Appendix \ref{A.I} that
\begin{eqnarray}
\sum_{m=k}^{n-1}(-1)^{n+m}\frac{2m+1}{(n-m)(n+m+1)}
&=& -\psi(2n+1)+\psi(n+k+1)+\psi(n+1)-\psi(n-k+1)
\nonumber \\
&& -\,\psi\left(\left\lfloor\frac{n+k}{2}\right\rfloor+1\right)
+\psi\left(\left\lfloor\frac{n-k}{2}\right\rfloor+1\right)
\nonumber \\
&& \hspace*{10em} (0\leqslant k\leqslant n-1),
\label{2.38}
\end{eqnarray}
where $\lfloor x\rfloor=\max\{n\in\mathbb{Z}:\:n\leqslant x\}$. Use
of Eq.\ (\ref{2.38}) casts Eq.\ (\ref{2.37}) into the final form
\begin{eqnarray}
c_{nk} &=& (-1)^{n+k}\frac{8(2k+1)}{(n-k)(n+k+1)}
\nonumber \\
&& \quad \times\bigg[\psi(n+k+1)-\psi(n-k+1)
-\psi\left(\left\lfloor\frac{n+k}{2}\right\rfloor+1\right)
+\psi\left(\left\lfloor\frac{n-k}{2}\right\rfloor+1\right)\bigg]
\nonumber \\
&& -\,\frac{4(2k+1)^{2}}{(n-k)^{2}(n+k+1)^{2}}
-(-1)^{n+k}\frac{4(2n+1)(2k+1)}{(n-k)^{2}(n+k+1)^{2}}
\qquad (0\leqslant k\leqslant n-1).
\nonumber \\
&&
\label{2.39}
\end{eqnarray}

Equation (\ref{2.39}) says nothing about the coefficient $c_{nn}$.
But from Eqs.\ (\ref{2.23}), (\ref{2.21}) and (\ref{2.33}) it can be
deduced that $c_{nn}$ may be expressed as
\begin{equation}
c_{nn}=-\sum_{k=0}^{n-1}c_{nk}.
\label{2.40}
\end{equation}
This implies that the polynomial $C_{n}(z)$ may be written as
\begin{equation}
C_{n}(z)=\sum_{k=0}^{n-1}c_{nk}[P_{k}(z)-P_{n}(z)],
\label{2.41}
\end{equation}
or explicitly, if the result in Eq.\ (\ref{2.39}) is used, as
\begin{eqnarray}
C_{n}(z) &=& 4\sum_{k=0}^{n-1}(-1)^{n+k}
\frac{2k+1}{(n-k)(n+k+1)}
\nonumber \\
&& \times\,\bigg\{2\bigg[\psi(n+k+1)-\psi(n-k+1)
-\psi\left(\left\lfloor\frac{n+k}{2}\right\rfloor+1\right)
+\psi\left(\left\lfloor\frac{n-k}{2}\right\rfloor+1\right)\bigg]
\nonumber \\
&& \quad -\,(-1)^{n+k}\frac{2k+1}{(n-k)(n+k+1)}
-\frac{2n+1}{(n-k)(n+k+1)}\bigg\}[P_{k}(z)-P_{n}(z)].
\label{2.42}
\end{eqnarray}

A bit different formula for $C_{n}(z)$ is obtained if the coefficient
$c_{nn}$ is expressed in a closed form. To find the latter, we
combine Eqs.\ (\ref{2.37}) and (\ref{2.40}) and write
\begin{eqnarray}
c_{nn} &=& -8[\psi(2n+1)-\psi(n+1)]
\sum_{k=0}^{n-1}(-1)^{n+k}\frac{2k+1}{(n-k)(n+k+1)}
\nonumber \\
&& -\,8\sum_{k=0}^{n-1}(-1)^{k}\frac{2k+1}{(n-k)(n+k+1)}
\sum_{m=k}^{n-1}(-1)^{m}\frac{2m+1}{(n-m)(n+m+1)}
\nonumber \\
&& +\,4\sum_{k=0}^{n-1}\frac{(2k+1)^{2}}{(n-k)^{2}(n+k+1)^{2}}
+4\sum_{k=0}^{n-1}(-1)^{n+k}
\frac{(2n+1)(2k+1)}{(n-k)^{2}(n+k+1)^{2}}.
\label{2.43}
\end{eqnarray}
The sums appearing in Eq.\ (\ref{2.43}) are evaluated in individual
subsections of the appendix, where it is found that
\begin{eqnarray}
&& \sum_{k=0}^{n-1}(-1)^{n+k}\frac{2k+1}{(n-k)(n+k+1)}
=-\psi(2n+1)+\psi(n+1),
\label{2.44}
\end{eqnarray}
\begin{eqnarray}
&& \hspace*{-5em}
\sum_{k=0}^{n-1}(-1)^{k}\frac{2k+1}{(n-k)(n+k+1)}
\sum_{m=k}^{n-1}(-1)^{m}\frac{2m+1}{(n-m)(n+m+1)}
\nonumber \\
&& =\,\frac{\pi^{2}}{12}-\frac{\gamma}{2n+1}
+\frac{1}{2}[\psi(2n+1)-\psi(n+1)]^{2}
-\frac{1}{2n+1}\psi(2n+1)-\frac{1}{2}\psi_{1}(2n+1),
\nonumber \\
&&
\label{2.45}
\end{eqnarray}
\begin{equation}
\sum_{k=0}^{n-1}\frac{(2k+1)^{2}}{(n-k)^{2}(n+k+1)^{2}}
=\frac{\pi^{2}}{6}-\frac{2\gamma}{2n+1}-\frac{2}{2n+1}\psi(2n+1)
-\psi_{1}(2n+1)
\label{2.46}
\end{equation}
and
\begin{equation}
\sum_{k=0}^{n-1}(-1)^{n+k}\frac{(2n+1)(2k+1)}{(n-k)^{2}(n+k+1)^{2}}
=-\frac{\pi^{2}}{12}+\psi_{1}(2n+1)-\frac{1}{2}\psi_{1}(n+1),
\label{2.47}
\end{equation}
with $\gamma$ standing for the Euler--Mascheroni constant and with
$\psi_{1}(\zeta)=\mathrm{d}\psi(\zeta)/\mathrm{d}\zeta$ being the
trigamma function. Plugging the results (\ref{2.44})--(\ref{2.47})
into the right-hand side of Eq.\ (\ref{2.43}) furnishes the
coefficient $c_{nn}$ in the compact form
\begin{equation}
c_{nn}=-\frac{\pi^{2}}{3}+4[\psi(2n+1)-\psi(n+1)]^{2}
+4\psi_{1}(2n+1)-2\psi_{1}(n+1).
\label{2.48}
\end{equation}
Hence, by virtue of Eqs.\ (\ref{2.33}), (\ref{2.39}) and
(\ref{2.48}), we eventually arrive at the sought formula
\begin{eqnarray}
C_{n}(z) &=& \left\{-\frac{\pi^{2}}{3}+4[\psi(2n+1)-\psi(n+1)]^{2}
+4\psi_{1}(2n+1)-2\psi_{1}(n+1)\right\}P_{n}(z)
\nonumber \\
&& +\,4\sum_{k=0}^{n-1}(-1)^{n+k}\frac{2k+1}{(n-k)(n+k+1)}
\bigg\{2\bigg[\psi(n+k+1)-\psi(n-k+1)
\nonumber \\
&& \qquad -\,\psi\left(\left\lfloor\frac{n+k}{2}\right\rfloor+1\right)
+\psi\left(\left\lfloor\frac{n-k}{2}\right\rfloor+1\right)\bigg]
\nonumber \\
&& \quad -\,(-1)^{n+k}\frac{2k+1}{(n-k)(n+k+1)}
-\frac{2n+1}{(n-k)(n+k+1)}\bigg\}P_{k}(z),
\label{2.49}
\end{eqnarray}
alternative to the one in Eq.\ (\ref{2.42}).
%
%
\subsection{Explicit expressions for 
$[\partial^{2}P_{\nu}(z)/\partial\nu^{2}]_{\nu=n}$ with $0\leqslant
n\leqslant3$}
\label{II.5}
It may be of interest to see how the derivatives
$[\partial^{2}P_{\nu}(z)/\partial\nu^{2}]_{\nu=n}$ look explicitly
for several lowest values of $n$. From Eqs.\ (\ref{2.12}),
(\ref{2.32}) and (\ref{2.49}), for $0\leqslant n\leqslant3$ we find
that
\begin{subequations}
\begin{eqnarray}
\frac{\partial^{2}P_{\nu}(z)}{\partial\nu^{2}}\bigg|_{\nu=0}
&=& -2\Li_{2}\frac{1-z}{2}, 
\label{2.50a}
\\*[1ex]
\frac{\partial^{2}P_{\nu}(z)}{\partial\nu^{2}}\bigg|_{\nu=1}
&=& -2z\Li_{2}\frac{1-z}{2}+(2z+2)\ln\frac{z+1}{2}-2z+2,
\label{2.50b}
\\*[1ex]
\frac{\partial^{2}P_{\nu}(z)}{\partial\nu^{2}}\bigg|_{\nu=2}
&=& (-3z^{2}+1)\Li_{2}\frac{1-z}{2}
+\left(\frac{7}{2}z^{2}+3z-\frac{1}{2}\right)\ln\frac{z+1}{2}
-\frac{11}{4}z^{2}+\frac{5}{2}z+\frac{1}{4},
\nonumber \\
&&
\label{2.50c}
\\*[1ex]
\frac{\partial^{2}P_{\nu}(z)}{\partial\nu^{2}}\bigg|_{\nu=3}
&=& (-5z^{3}+3z)\Li_{2}\frac{1-z}{2}
+\left(\frac{37}{6}z^{3}+5z^{2}-\frac{5}{2}z
-\frac{4}{3}\right)\ln\frac{z+1}{2}
\nonumber \\
&& -\,\frac{155}{36}z^{3}+\frac{23}{6}z^{2}
+\frac{19}{12}z-\frac{10}{9}.
\label{2.50d}
\end{eqnarray}
\label{2.50}%
\end{subequations}
%
%
\section{The derivatives $[\partial Q_{\nu}(z)/\partial\nu]_{\nu=n}$
and $[\partial Q_{\nu}(x)/\partial\nu]_{\nu=n}$}
\label{III}
\setcounter{equation}{0}
In Refs.\ \cite{Szmy06,Szmy07}, we exploited representations of the
first-order derivatives $[\partial P_{\nu}(z)/\partial\nu]_{\nu=n}$
found therein to obtain expressions for the Legendre functions of the
second kind $Q_{n}(z)$, with $n\in\mathbb{N}_{0}$, both for
$z\in\mathbb{C}\setminus[-1,1]$ and for $z=x\in(-1,1)$. Below we
shall show that the knowledge of the second-order derivatives
$[\partial^{2}P_{\nu}(z)/\partial\nu^{2}]_{\nu=n}$ allows one to
obtain explicit formulas for the first-order derivatives $[\partial
Q_{\nu}(z)/\partial\nu]_{\nu=n}$ and $[\partial
Q_{\nu}(x)/\partial\nu]_{\nu=n}$, again with $n\in\mathbb{N}_{0}$.
\subsection{The derivatives $[\partial
Q_{\nu}(z)/\partial\nu]_{\nu=n}$ for $z\in\mathbb{C}\setminus[-1,1]$}
\label{III.1}
\subsubsection{The general form of $[\partial
Q_{\nu}(z)/\partial\nu]_{\nu=n}$}
\label{III.1.1}
The Legendre function of the second kind, $Q_{\nu}(z)$, may be
defined in terms of the function of the first kind through the
formula
\begin{equation}
Q_{\nu}(z)=\frac{\pi}{2}
\frac{\mathrm{e}^{\mp\mathrm{i}\pi\nu}P_{\nu}(z)-P_{\nu}(-z)}
{\sin(\pi\nu)}
\qquad (\Imag z\gtrless0).
\label{3.1}
\end{equation}
Hence, it follows that
\begin{eqnarray}
&& \frac{\partial
Q_{\nu}(z)}{\partial\nu}=\frac{\pi}{2\sin^{2}(\pi\nu)}
\left\{-\pi[P_{\nu}(z)-P_{\nu}(-z)\cos(\pi\nu)]
+\left[\mathrm{e}^{\mp\mathrm{i}\pi\nu}
\frac{\partial P_{\nu}(z)}{\partial\nu}
-\frac{\partial P_{\nu}(-z)}{\partial\nu}\right]\sin(\pi\nu)\right\}
\nonumber \\
&& \hspace*{30em} (\Imag z\gtrless0).
\label{3.2}
\end{eqnarray}
From this, for $\nu=n\in\mathbb{N}_{0}$, with the use of the
L'Hospital rule, we obtain
\begin{equation}
\frac{\partial Q_{\nu}(z)}{\partial\nu}\bigg|_{\nu=n}
=-\frac{\pi^{2}}{4}P_{n}(z)\mp\frac{\mathrm{i}\pi}{2}
\frac{\partial P_{\nu}(z)}{\partial\nu}\bigg|_{\nu=n}
+\frac{1}{4}\frac{\partial^{2}P_{\nu}(z)}
{\partial\nu^{2}}\bigg|_{\nu=n}
-\frac{(-1)^{n}}{4}\frac{\partial^{2}P_{\nu}(-z)}{\partial\nu^{2}}
\bigg|_{\nu=n}
\qquad (\Imag z\gtrless0).
\label{3.3}
\end{equation}
If in the above formula the second-order derivatives
$[\partial^{2}P_{\nu}(\pm z)/\partial\nu^{2}]_{\nu=n}$ are
substituted with expressions following from Eq.\ (\ref{2.12}) and the
first-order derivative $[\partial P_{\nu}(z)/\partial\nu]_{\nu=n}$ is
replaced by the right-hand side of Eq.\ (\ref{1.1}), this yields
$[\partial Q_{\nu}(z)/\partial\nu]_{\nu=n}$ in the form
\begin{eqnarray}
\frac{\partial Q_{\nu}(z)}{\partial\nu}\bigg|_{\nu=n}
&=& \frac{1}{2}P_{n}(z)
\left(\Li_{2}\frac{z+1}{2}-\Li_{2}\frac{1-z}{2}\right)
+\left[\frac{1}{4}B_{n}(z)\mp\frac{\mathrm{i}\pi}{2}P_{n}(z)\right]
\ln\frac{z+1}{2}
\nonumber \\
&& -\,\frac{(-1)^{n}}{4}B_{n}(-z)\ln\frac{1-z}{2}
-\frac{\pi^{2}}{4}P_{n}(z)
\mp\frac{\mathrm{i}\pi}{2}R_{n}(z)+\frac{1}{4}C_{n}(z)
-\frac{(-1)^{n}}{4}C_{n}(-z)
\nonumber \\
&& \hspace*{25em} (\Imag z\gtrless0).
\label{3.4}
\end{eqnarray}
A more elegant expression for $[\partial
Q_{\nu}(z)/\partial\nu]_{\nu=n}$ follows if the dilogarithm
$\Li_{2}[(z+1)/2]$ is eliminated from Eq.\ (\ref{3.4}) with the aid
of the Euler's identity \cite[Eq.\ (1.11)]{Lewi81}
\begin{equation}
\Li_{2}z+\Li_{2}(1-z)=\frac{\pi^{2}}{6}-\ln z\ln(1-z),
\label{3.5}
\end{equation}
the relation
\begin{equation}
1-z=\mathrm{e}^{\mp\mathrm{i}\pi}(z-1)
\qquad (\Imag z\gtrless0)
\label{3.6}
\end{equation}
and the result in Eq.\ (\ref{2.31}). Proceeding in that way, one
eventually finds that
\begin{eqnarray}
\frac{\partial Q_{\nu}(z)}{\partial\nu}\bigg|_{\nu=n}
&=& -P_{n}(z)\Li_{2}\frac{1-z}{2}
-\frac{1}{2}P_{n}(z)\ln\frac{z+1}{2}\ln\frac{z-1}{2}
+\frac{1}{4}B_{n}(z)\ln\frac{z+1}{2}
\nonumber \\
&& -\,\frac{(-1)^{n}}{4}B_{n}(-z)\ln\frac{z-1}{2}
-\frac{\pi^{2}}{6}P_{n}(z)
+\frac{1}{4}C_{n}(z)-\frac{(-1)^{n}}{4}C_{n}(-z)
\qquad (n\in\mathbb{N}_{0}).
\nonumber \\
&&
\label{3.7}
\end{eqnarray}
\subsubsection{Explicit expressions for $[\partial
Q_{\nu}(z)/\partial\nu]_{\nu=n}$ with $0\leqslant n\leqslant3$}
\label{III.1.2}
Explicit forms of the derivatives $[\partial
Q_{\nu}(z)/\partial\nu]_{\nu=n}$ with $0\leqslant n\leqslant3$,
obtained from Eq.\ (\ref{3.7}) with the use of Eqs.\ (\ref{2.32}) and
(\ref{2.49}), are
\begin{subequations}
\begin{eqnarray}
\frac{\partial Q_{\nu}(z)}{\partial\nu}\bigg|_{\nu=0}
&=& -\Li_{2}\frac{1-z}{2}-\frac{1}{2}\ln\frac{z+1}{2}\ln\frac{z-1}{2}
-\frac{\pi^{2}}{6},
\label{3.8a}
\\*[1ex]
\frac{\partial Q_{\nu}(z)}{\partial\nu}\bigg|_{\nu=1}
&=& -z\Li_{2}\frac{1-z}{2}
-\frac{1}{2}z\ln\frac{z+1}{2}\ln\frac{z-1}{2}
\nonumber \\
&& +\,\left(\frac{1}{2}z+\frac{1}{2}\right)\ln\frac{z+1}{2}
+\left(-\frac{1}{2}z+\frac{1}{2}\right)\ln\frac{z-1}{2}
-\frac{\pi^{2}}{6}z+1,
\label{3.8b}
\\*[1ex]
\frac{\partial Q_{\nu}(z)}{\partial\nu}\bigg|_{\nu=2}
&=& \left(-\frac{3}{2}z^{2}+\frac{1}{2}\right)\Li_{2}\frac{1-z}{2}
+\left(-\frac{3}{4}z^{2}+\frac{1}{4}\right)
\ln\frac{z+1}{2}\ln\frac{z-1}{2}
\nonumber \\
&& +\,\left(\frac{7}{8}z^{2}+\frac{3}{4}z-\frac{1}{8}\right)
\ln\frac{z+1}{2}
+\left(-\frac{7}{8}z^{2}+\frac{3}{4}z+\frac{1}{8}\right)
\ln\frac{z-1}{2}-\frac{\pi^{2}}{4}z^{2}
+\frac{5}{4}z+\frac{\pi^{2}}{12},
\nonumber \\
&&
\label{3.8c}
\\*[1ex]
\frac{\partial Q_{\nu}(z)}{\partial\nu}\bigg|_{\nu=3}
&=& \left(-\frac{5}{2}z^{3}+\frac{3}{2}z\right)\Li_{2}\frac{1-z}{2}
+\left(-\frac{5}{4}z^{3}+\frac{3}{4}z\right)
\ln\frac{z+1}{2}\ln\frac{z-1}{2}
\nonumber \\
&& +\,\left(\frac{37}{24}z^{3}+\frac{5}{4}z^{2}
-\frac{5}{8}z-\frac{1}{3}\right)\ln\frac{z+1}{2}
+\left(-\frac{37}{24}z^{3}+\frac{5}{4}z^{2}
+\frac{5}{8}z-\frac{1}{3}\right)\ln\frac{z-1}{2}
\nonumber \\
&& -\,\frac{5\pi^{2}}{12}z^{3}+\frac{23}{12}z^{2}
+\frac{\pi^{2}}{4}z-\frac{5}{9}.
\label{3.8d}
\end{eqnarray}
\label{3.8}%
\end{subequations}
We find it remarkable that coefficients in the polynomial part of
$[\partial Q_{\nu}(z)/\partial\nu]_{\nu=n}$ are alternately
irrational and rational.
\subsection{The derivatives $[\partial
Q_{\nu}(x)/\partial\nu]_{\nu=n}$ for $-1<x<1$}
\label{III.2}
On the real interval $-1<x<1$, the Legendre function of the second
kind, $Q_{\nu}(x)$, is defined as the average of the limits
$Q_{\nu}(x+\mathrm{i}0)$ and $Q_{\nu}(x-\mathrm{i}0)$ resulting when
$z$ approaches $x$ from the upper ($\Imag z>0$) and lower ($\Imag
z<0$) half-planes, respectively. One has
\begin{equation}
Q_{\nu}(x)=\frac{1}{2}[Q_{\nu}(x+\mathrm{i}0)+Q_{\nu}(x-\mathrm{i}0)]
\qquad (-1<x<1),
\label{3.9}
\end{equation}
and consequently
\begin{equation}
\frac{\partial Q_{\nu}(x)}{\partial\nu}\bigg|_{\nu=n}
=\frac{1}{2}\frac{\partial Q_{\nu}(x+\mathrm{i}0)}
{\partial\nu}\bigg|_{\nu=n}
+\frac{1}{2}\frac{\partial Q_{\nu}(x-\mathrm{i}0)}
{\partial\nu}\bigg|_{\nu=n}.
\label{3.10}
\end{equation}
From this, with the use of Eq.\ (\ref{3.7}) and the identity
\begin{equation}
x-1\pm\mathrm{i}0=\mathrm{e}^{\pm\mathrm{i}\pi}(1-x)
\qquad (-1<x<1),
\label{3.11}
\end{equation}
one finds that
\begin{eqnarray}
\frac{\partial Q_{\nu}(x)}{\partial\nu}\bigg|_{\nu=n}
&=& -P_{n}(x)\Li_{2}\frac{1-x}{2}
-\frac{1}{2}P_{n}(x)\ln\frac{1+x}{2}\ln\frac{1-x}{2}
+\frac{1}{4}B_{n}(x)\ln\frac{1+x}{2}
\nonumber \\
&& -\,\frac{(-1)^{n}}{4}B_{n}(-x)\ln\frac{1-x}{2}
-\frac{\pi^{2}}{6}P_{n}(x)
+\frac{1}{4}C_{n}(x)-\frac{(-1)^{n}}{4}C_{n}(-x)
\qquad (n\in\mathbb{N}_{0}).
\nonumber \\
&&
\label{3.12}
\end{eqnarray}

There is no need to provide here explicit representations for the
derivatives $[\partial Q_{\nu}(x)/\partial\nu]_{\nu=n}$ for several
lowest non-negative values of $n$. As it is seen from Eqs.\
(\ref{3.7}) and (\ref{3.12}), such representations for $0\leqslant
n\leqslant3$ may be immediately deduced from Eqs.\
(\ref{3.8a})--(\ref{3.8d}) after the replacement of $z$ with $x$ is
made everywhere in the latter set of equations, except for the
logarithm $\ln[(z-1)/2]$, which is to be substituted with
$\ln[(1-x)/2]$.
\section*{Acknowledgments}
I wish to thank Dr.\ George P.\ Schramkowski for kindly communicating
to me the formula in Eq.\ (\ref{1.3}) and for the subsequent
inspiring correspondence.
%
%
\appendix
\renewcommand{\theequation}{\thesubsection.\arabic{equation}}
\section{Appendix: Proofs of summation formulas used in Sec.\
\ref{II.4}} 
\label{A} 
\setcounter{equation}{0}
\subsection{The summation formulas (\ref{2.38}) and (\ref{2.44})}
\label{A.I} 
\setcounter{equation}{0}
We denote
\begin{equation}
S_{1}=\sum_{k=m}^{n-1}(-1)^{n+k}\frac{2k+1}{(n-k)(n+k+1)}
\qquad (0\leqslant m\leqslant n-1).
\label{A.1}
\end{equation}
We have
\begin{equation}
S_{1}=\sum_{k=m}^{n-1}\frac{(-1)^{n+k}}{n-k}
-\sum_{k=m}^{n-1}\frac{(-1)^{n+k}}{n+k+1}
\label{A.2}
\end{equation}
and further
\begin{equation}
S_{1}=\sum_{k=1}^{n-m}\frac{(-1)^{k}}{k}
+\sum_{k=n+m+1}^{2n}\frac{(-1)^{k}}{k}
=\sum_{k=1}^{n-m}\frac{(-1)^{k}}{k}
+\sum_{k=1}^{2n}\frac{(-1)^{k}}{k}
-\sum_{k=1}^{n+m}\frac{(-1)^{k}}{k}.
\label{A.3}
\end{equation}
However, it is easy to show that
\begin{equation}
\sum_{k=1}^{N}\frac{(-1)^{k}}{k}
=-\sum_{k=1}^{N}\frac{1}{k}
+\sum_{k=1}^{\lfloor N/2\rfloor}\frac{1}{k}
\qquad (N\in\mathbb{N}_{0}),
\label{A.4}
\end{equation}
where $\lfloor x\rfloor=\max\{n\in\mathbb{Z}:\:n\leqslant x\}$ stands
for the integer part of $x$. Since
\begin{equation}
\sum_{k=1}^{N}\frac{1}{k}
=\psi(N+1)-\psi(1)
\qquad (N\in\mathbb{N}_{0}),
\label{A.5}
\end{equation}
with $\psi(z)=\mathrm{d}\ln\Gamma(z)/\mathrm{d}z$ being the digamma
function, Eq.\ (\ref{A.4}) may be rewritten in the form
\begin{equation}
\sum_{k=1}^{N}\frac{(-1)^{k}}{k}=-\psi(N+1)
+\psi\left(\bigg\lfloor\frac{N}{2}\bigg\rfloor+1\right)
\qquad (N\in\mathbb{N}_{0}).
\label{A.6}
\end{equation}
Application of the result (\ref{A.6}) to each of the three sums on
the extreme right-hand side of Eq.\ (\ref{A.3}) gives finally
\begin{eqnarray}
\sum_{k=m}^{n-1}(-1)^{n+k}\frac{2k+1}{(n-k)(n+k+1)}
&=& -\psi(2n+1)+\psi(n+m+1)+\psi(n+1)-\psi(n-m+1)
\nonumber \\
&& -\,\psi\left(\bigg\lfloor\frac{n+m}{2}\bigg\rfloor+1\right)
+\psi\left(\bigg\lfloor\frac{n-m}{2}\bigg\rfloor+1\right)
\nonumber \\
&& \hspace*{13em} (0\leqslant m\leqslant n-1).
\label{A.7}
\end{eqnarray}
After $k$ is interchanged with $m$, Eq.\ (\ref{A.7}) becomes
identical with Eq.\ (\ref{2.38}).

For $m=0$, Eq.\ (\ref{A.7}) becomes
\begin{equation}
\sum_{k=0}^{n-1}(-1)^{n+k}\frac{2k+1}{(n-k)(n+k+1)}
=-\psi(2n+1)+\psi(n+1),
\label{A.8}
\end{equation}
which is Eq.\ (\ref{2.44}).
%
%
\subsection{The summation formula (\ref{2.45})}
\label{A.II}
\setcounter{equation}{0}
We denote
\begin{equation}
S_{2}=\sum_{k=0}^{n-1}(-1)^{k}\frac{2k+1}{(n-k)(n+k+1)}
\sum_{m=k}^{n-1}(-1)^{m}\frac{2m+1}{(n-m)(n+m+1)}.
\label{A.9}
\end{equation}
Application of the identity
\begin{equation}
\sum_{k=N_{1}}^{N_{2}}f_{k}\sum_{m=k}^{N_{2}}f_{m}
=\frac{1}{2}\left(\sum_{k=N_{1}}^{N_{2}}f_{k}\right)^{2}
+\frac{1}{2}\sum_{k=N_{1}}^{N_{2}}f_{k}^{2}
\qquad (N_{1}\leqslant N_{2})
\label{A.10}
\end{equation}
transforms Eq.\ (\ref{A.9}) into
\begin{equation}
S_{2}=\frac{1}{2}\left[\sum_{k=0}^{n-1}(-1)^{k}
\frac{2k+1}{(n-k)(n+k+1)}\right]^{2}
+\frac{1}{2}\sum_{k=0}^{n-1}\frac{(2k+1)^{2}}{(n-k)^{2}(n+k+1)^{2}},
\label{A.11}
\end{equation}
from which, with the help of Eqs.\ (\ref{A.8}) and (\ref{A.22}),
we obtain
\begin{eqnarray}
&& \hspace*{-5em}
\sum_{k=0}^{n-1}(-1)^{k}\frac{2k+1}{(n-k)(n+k+1)}
\sum_{m=k}^{n-1}(-1)^{m}\frac{2m+1}{(n-m)(n+m+1)}
\nonumber \\
&& =\,\frac{\pi^{2}}{12}-\frac{\gamma}{2n+1}
+\frac{1}{2}[\psi(2n+1)-\psi(n+1)]^{2}
-\frac{1}{2n+1}\psi(2n+1)-\frac{1}{2}\psi_{1}(2n+1),
\nonumber \\
&&
\label{A.12}
\end{eqnarray}
which is Eq.\ (\ref{2.45}).

To prove the identity (\ref{A.10}), we write the obvious chain of
equalities ($N_{1}\leqslant N_{2}$ is assumed)
\begin{equation}
\left(\sum_{k=N_{1}}^{N_{2}}f_{k}\right)^{2}
=\sum_{k=N_{1}}^{N_{2}}f_{k}\sum_{m=N_{1}}^{N_{2}}f_{m}
=\sum_{k=N_{1}}^{N_{2}}f_{k}\sum_{m=N_{1}}^{k}f_{m}
+\sum_{k=N_{1}}^{N_{2}}f_{k}\sum_{m=k}^{N_{2}}f_{m}
-\sum_{k=N_{1}}^{N_{2}}f_{k}^{2}.
\label{A.13}
\end{equation}
Manipulating with the first term on the extreme right-hand side of
Eq.\ (\ref{A.13}), we have
\begin{equation}
\sum_{k=N_{1}}^{N_{2}}f_{k}\sum_{m=N_{1}}^{k}f_{m}
=\sum_{m=N_{1}}^{N_{2}}f_{m}\sum_{k=m}^{N_{2}}f_{k}
=\sum_{k=N_{1}}^{N_{2}}f_{k}\sum_{m=k}^{N_{2}}f_{m}.
\label{A.14}
\end{equation}
Plugging the result (\ref{A.14}) into Eq.\ (\ref{A.13}), we obtain
\begin{equation}
\left(\sum_{k=N_{1}}^{N_{2}}f_{k}\right)^{2}
=2\sum_{k=N_{1}}^{N_{2}}f_{k}\sum_{m=k}^{N_{2}}f_{m}
-\sum_{k=N_{1}}^{N_{2}}f_{k}^{2},
\label{A.15}
\end{equation}
from which the identity in Eq.\ (\ref{A.10}) follows immediately.
%
%
\subsection{The summation formula (\ref{2.46})}
\label{A.III}
\setcounter{equation}{0}
We denote
\begin{equation}
S_{3}=\sum_{k=0}^{n-1}\frac{(2k+1)^{2}}{(n-k)^{2}(n+k+1)^{2}}.
\label{A.16}
\end{equation}
If we carry out the partial fraction decomposition of the summand, we
have
\begin{equation}
S_{3}=\sum_{k=0}^{n-1}\frac{1}{(n-k)^{2}}
+\sum_{k=0}^{n-1}\frac{1}{(n+k+1)^{2}}
-\frac{2}{2n+1}\sum_{k=0}^{n-1}\frac{1}{n-k}
-\frac{2}{2n+1}\sum_{k=0}^{n-1}\frac{1}{n+k+1}
\label{A.17}
\end{equation}
and further, after obvious rearrangements,
\begin{equation}
S_{3}=\sum_{k=0}^{2n-1}\frac{1}{(k+1)^{2}}
-\frac{2}{2n+1}\sum_{k=1}^{2n}\frac{1}{k}.
\label{A.18}
\end{equation}
Now, it holds that
\begin{equation}
\sum_{k=0}^{N-1}\frac{1}{(k+1)^{2}}
=\sum_{k=0}^{\infty}\frac{1}{(k+1)^{2}}
-\sum_{k=0}^{\infty}\frac{1}{(k+N+1)^{2}}
=\psi_{1}(1)-\psi_{1}(N+1)
\qquad (N\in\mathbb{N}_{0}),
\label{A.19}
\end{equation}
where $\psi_{1}(z)=\mathrm{d}\psi(z)/\mathrm{d}z$ is the trigamma
function. On employing Eqs.\ (\ref{A.19}) and (\ref{A.5}) in Eq.\
(\ref{A.18}), after using the well-known relations
\begin{equation}
\psi(1)=-\gamma
\label{A.20}
\end{equation}
(here and below $\gamma$ stands for the Euler--Mascheroni constant)
and
\begin{equation}
\psi_{1}(1)=\frac{\pi^{2}}{6},
\label{A.21}
\end{equation}
we finally obtain
\begin{equation}
\sum_{k=0}^{n-1}\frac{(2k+1)^{2}}{(n-k)^{2}(n+k+1)^{2}}
=\frac{\pi^{2}}{6}-\frac{2\gamma}{2n+1}
-\frac{2}{2n+1}\psi(2n+1)-\psi_{1}(2n+1),
\label{A.22}
\end{equation}
which is Eq.\ (\ref{2.46}).
%
%
\subsection{The summation formula (\ref{2.47})}
\label{A.IV}
\setcounter{equation}{0}
We denote
\begin{equation}
S_{4}=\sum_{k=0}^{n-1}(-1)^{n+k}\frac{(2n+1)(2k+1)}
{(n-k)^{2}(n+k+1)^{2}}.
\label{A.23}
\end{equation}
A partial-fraction decomposition of the summand gives
\begin{equation}
S_{4}=\sum_{k=0}^{n-1}\frac{(-1)^{n+k}}{(n-k)^{2}}
-\sum_{k=0}^{n-1}\frac{(-1)^{n+k}}{(n+k+1)^{2}}.
\label{A.24}
\end{equation}
With a little bit of algebra on the right-hand side of Eq.\
(\ref{A.24}), we obtain
\begin{equation}
S_{4}=\sum_{k=0}^{2n-1}\frac{(-1)^{k+1}}{(k+1)^{2}}
\label{A.25}
\end{equation}
and further
\begin{equation}
S_{4}=-\sum_{k=0}^{2n-1}\frac{1}{(k+1)^{2}}
+\frac{1}{2}\sum_{k=0}^{n-1}\frac{1}{(k+1)^{2}}.
\label{A.26}
\end{equation}
From this, with reference to Eqs.\ (\ref{A.19}) and (\ref{A.21}),
we eventually arrive at
\begin{equation}
\sum_{k=0}^{n-1}(-1)^{n+k}\frac{(2n+1)(2k+1)}{(n-k)^{2}(n+k+1)^{2}}
=-\frac{\pi^{2}}{12}+\psi_{1}(2n+1)-\frac{1}{2}\psi_{1}(n+1),
\label{A.27}
\end{equation}
which is Eq.\ (\ref{2.47}).
%
%

%
\end{document}